\documentclass[12pt, a4paper]{amsart}
\usepackage{hyperref, fullpage, amsmath, amsfonts, amsthm, amssymb, stmaryrd}

\theoremstyle{plain}
\newtheorem{Theorem}{Theorem}[section]

\theoremstyle{definition}

\theoremstyle{remark}

\author{Arno Fehm}
\address{School of Mathematics, Alan Turing Building, The University of Manchester, Oxford Road, Manchester, M13 9PL, United Kingdom}
\email{arno.fehm@manchester.ac.uk}

\author{Franziska Jahnke}
\address{Institut f\"ur Mathematische Logik und Grundlagenforschung,
University of M\"unster, Einsteinstr.\;62, 48149 M\"unster, Germany}
\email{franziska.jahnke@wwu.de} 
\thanks{The second author was partially supported by SFB 878.}


\begin{document}

\title{Recent Progress on Definability of Henselian Valuations}

\begin{abstract}
Although the study of the definability of henselian valuations has a long history starting with J. Robinson,
most of the results in this area were proven during the last few years. 
We survey these results which address 
the definability of concrete henselian valuations,
the existence of definable henselian valuations on a given field, 
and questions of uniformity and quantifier complexity.
\end{abstract}

\maketitle


This is a survey about definable henselian valuations.
The main motivation to study definable henselian valuations has always come from the desire to understand definable subsets of a field,
in particular for applications regarding decidability. 
This subject has recently received a lot of attention and numerous results have been proven (including by the authors) 
in the last five years.
One reason for this increased activity is the Shelah-Hasson conjecture
that every NIP field which is neither separably closed nor real closed admits a nontrivial definable henselian valuation.
Although there are a number of very interesting and deep results on definable valuations which are not henselian,
in particular valuations on global fields, 
we restrict ourselves to results on henselian valuations.
Similarly, for reasons of limited space, while we do try to mention every paper in this area known to us,
we had to make a selection of the results that we present here.

We start by recalling some basic definitions:
For a valued field $(K,v)$ we denote the valuation ring by $\mathcal{O}_v$, 
the value group by $vK$
and the residue field by $Kv$.
The valuation $v$ is {\bf henselian} if it extends uniquely to every finite extension of $K$,
and {\bf $p$-henselian}, for $p$ a prime number, if it extends uniquely to every finite Galois extension of degree $p$;
for more details and equivalent statements see \cite[Theorem 4.1.3]{EP}, respectively \cite[Theorem 4.2.2, 4.2.3]{EP}.
Given two valuations $v$ and $w$ on a field $K$, we call $v$ coarser than $w$ if $\mathcal{O}_w \subseteq \mathcal{O}_v$ holds.
The {\bf canonical henselian} (resp.~{\bf canonical $p$-henselian}) valuation on $K$
is the coarsest henselian (resp.~$p$-henselian) valuation on $K$ with separably closed (resp.~$p$-closed, i.e.~not admitting a Galois extension of degree $p$) residue field, if any such exists,
and otherwise the finest henselian (resp.~$p$-henselian) valuation on $K$.
The canonical henselian (resp.~$p$-henselian) valuation exists and is nontrivial whenever 
$K$ admits a nontrivial henselian (resp.~$p$-henselian) valuation and is not separably closed (resp.~not $p$-closed),
cf.~\cite[\S4.4]{EP} (resp.~\cite[Prop.~3.1]{Koe95}).
If $K$ admits some nontrivial henselian (resp.~$p$-henselian) valuation, 
we also call $K$ henselian (resp.~$p$-henselian).

We say that a valuation $v$ on $K$ is {\bf definable} if $\mathcal{O}_v$ is a
definable subset of $K$ in the language $\mathcal{L}_{\rm ring}$ of rings,
i.e., if there exists an $\mathcal{L}_{\rm ring}$-formula $\phi(x)$ such that for $a\in K$,
$$
 a\in\mathcal{O}_v \Longleftrightarrow K\models\phi(a),
$$
in which case we also say that $\phi$ defines $v$.
For a subset $A$ of $K$ we call $\phi(x)$ an $A$-formula if it uses parameters only from $A$,
and we call $\phi(x)$ an $\exists$-formula (resp.~$\exists\forall$-formula etc.) if
it is equivalent to a formula in prenex normal form with precisely these quantifiers.
We then say that $v$ is $A$-definable, $\exists$-$A$-definable etc.~if there exists a corresponding formula defining $v$.

This survey is organized as follows:
in Section \ref{sec:1} we give results regarding the definability of a given henselian valuation on a field $K$,
many of which are classical.
Section \ref{sec:2} treats the existence of {\em some} nontrivial definable henselian valuation on a given field,
including very recent results on fields with certain model-theoretic properties, in particular NIP.
Finally, in Section \ref{sec:3}, we 
survey a number of results from the last three years
discussing how the defining formulae can be optimized with respect to uniformity and quantifier complexity.

\section{Definability of a given henselian valuation} 
The first definability results of henselian valuations were shown 50 years ago. \label{sec:1}
The main aim was to reduce the decidability of some field to the decidability
of a subring or subfield. Julia Robinson showed the following theorem which implies
that the $\mathcal{L}_\textrm{ring}$-theory of $\mathbb{Q}_p$ is decidable if
and only if the $\mathcal{L}_\textrm{ring}$-theory of 
$\mathbb{Z}_p$ is decidable: 
\begin{Theorem}[J.\ Robinson, {\cite[\S2]{Robinson}}]\label{Robinson}
The $p$-adic valuation on $\mathbb{Q}_p$ is defined by the formula
$$
 \exists y\;(y^2=1+kx^2)
$$
where $k=p$ if $p$ is odd, and $k=8$ if $p=2$.
\end{Theorem}

Similarly, Ax gave a definition of the ring $F[[t]]$ in $F((t))$ to conclude
that, in the case ${\rm char}(F)=0$, 
the $\mathcal{L}_\textrm{ring}$-theory of $F((t))$ is decidable if and only if
the $\mathcal{L}_\textrm{ring}$-theory of $F$ is decidable.
\begin{Theorem}[Ax, \cite{Ax}]\label{Ax}
For any field $F$, the power series valuation on the field of Laurent series $F((t))$ 
is defined by
\begin{align*}
 \exists w, y \forall u, x_1, x_2 \exists z \forall y_1, y_2\, [&(z^m=1+wx_1^mx_2^m\vee y_1^m\neq 1+wx_1^m \vee y_2^m\neq 1+wx_2^m)
\\&\wedge u^m\neq w\wedge y^m=1+wx^m]
\end{align*}
where $m>1$ and ${\rm char}(F)\nmid m$.
\end{Theorem}

Moreover, Ax observed that both Julia Robinson's and his formulae work for arbitrary henselian valued fields with value group a $\mathbb{Z}$-group (i.e., when 
$vK \equiv \mathbb{Z}$ as ordered abelian groups).
A first generalization of these results was proven by Koenigsmann (\cite[Lemma 3.6]{Koe04}) who showed\footnote{Note that 
his proof relies on a result from \cite{Koe94} which does not quite work. 
However, the citation in question can be replaced by Theorem \ref{uniformp} below.} 
that any henselian valuation with an archimedean non-divisible value group is $\emptyset$-definable.  
In \cite[Corollary 2]{Hong}, Hong points out that Julia Robinson's formula can be generalized to henselian valued fields with discrete value group (i.e., when the value group has a least positive element).
Improving on this further and incorporating Ax's idea, Hong finally proves:
\begin{Theorem}[Hong, {\cite[Theorem 4]{Hong}}]
Let $(K,v)$ be a henselian valued field. \label{Hong}
If $vK$ is regular and non-divisible, then $v$ is $\emptyset$-definable.
\end{Theorem}
Here, an ordered abelian group is called \emph{regular} if every quotient by a non-zero convex subgroup is divisible.
In particular, all $\mathbb{Z}$-groups and all archimedean ordered abelian groups are regular.

All the above definitions of henselian valuations 
were obtained using properties of the value group of said valuation. 
On the other hand, one can also define a henselian valuation using specifics
about its residue field, as asserted by the next theorem.
\begin{Theorem}[Jahnke--Koenigsmann, {\cite[Proposition 3.1, Corollaries 3.3 and 3.8]{JK12}}] \label{JK12}
Let $(K,v)$ be a henselian valued field.
Then $v$ is $\emptyset$-definable if the residue field $Kv$ satisfies one of the following conditions:
\begin{enumerate}
\item for some prime $p>2$, $Kv$ is not $p$-henselian and not $p$-closed.
\item $Kv$ is hilbertian.
\item $Kv$ is pseudo-algebraically closed and not separably closed.
\end{enumerate}
\end{Theorem}

See \cite[Chapters 11 and 12]{FJ} for definitions of pseudo-algebraically closed and hilbertian fields.
In particular, condition (1) applies to the case where $Kv$ is {\em finite},
and condition (2) applies in particular to the case where $Kv$ is a {\em number field}.

As discussed in section \ref{sec:2} below, the canonical henselian valuation is not necessarily definable. The canonical $p$-henselian valuation 
is however usually definable. More precisely, 
correcting a mistake in \cite{Koe95}, Jahnke and Koenigsmann show:

\begin{Theorem}[Jahnke--Koenigsmann, {\cite[Theorem 3.1]{JK15b}}]\label{uniformp}
For every prime $p>2$ there exists a $\emptyset$-formula 
which defines the canonical $p$-henselian valuation 
on any field $K$ with either ${\rm char}(K)=p$ or $\mu_p\subseteq K$, where $\mu_p$ denotes the set of $p$th roots
of unity.
\end{Theorem}

Note that there is also a version for $p=2$, see \cite{JK15b} for more details.

\section{Definable henselian valuations on a given henselian field}

A field that admits a henselian valuation may in fact admit a multitude of them,
some of which might be definable, while others are not.
This section focuses on the question of when a henselian field admits any definable henselian valuation at all. \label{sec:2}

A first observation in this direction is that not every henselian field admits a nontrivial definable valuation:
Every nontrivial valuation on an algebraically closed, or, more generally, separably closed, field is henselian, although never definable\footnote{This follows from the stability of separably closed fields, see \cite{Wood}.}.
Also every non-archimedean real closed field carries nontrivial henselian valuations, none of which are definable.
There are, however, more henselian fields without definable henselian valuations:

\begin{Theorem}[Prestel--Ziegler, {\cite[p.~338]{PZ}}]\label{PZ}
There exists a henselian valued field $K$ of characteristic zero which is neither algebraically closed nor real closed
and does not admit any nontrivial $\emptyset$-definable henselian valuation.
\end{Theorem}

However, building on ideas from \cite{AEJ},
Koenigsmann sketched a proof in \cite{Koe94} that a henselian valued field which is neither separably closed nor real closed
always admits a definable valuation that at least induces the same topology as a henselian valuation.
There is no published article containing these 
results, and the preprint has some flaws.
However, it has been very influential, and has paved the way for a number of the current advances we discuss in this survey 
(in particular those in \cite{Hong}, \cite{JK15b}, \cite{JK16} and \cite{Krupinski}).
Koenigsmann's approach can be made rigorous by using results from \cite{JK15b}, see \cite{Dupont}.

Refinements of the Prestel--Ziegler construction have recently been given in \cite{FJ15}, \cite{JK16} and \cite{AJ}.
In particular, Theorem \ref{PZ} can in fact be strengthened from ``$\emptyset$-definable'' to ``definable'', see \cite[Example 6.2]{JK16}. 
The relations between admitting $\emptyset$-definable and definable henselian valuations have been explored further in \cite{AJ}.
They also show:

\begin{Theorem}[Anscombe--Jahnke, {\cite[Proposition 3.4]{AJ}}]
Let $K$ be a field which is not algebraically closed and whose canonical henselian valuation has residue characteristic zero.
If every field elementarily equivalent to $K$ is henselian, then $K$ admits a nontrivial $\emptyset$-definable henselian valuation.
\end{Theorem}

While in the previous section we listed 
results that give sufficient conditions on the residue field or the value group
of a henselian valuation 
to imply that this particular valuation is definable,
the following two results
give sufficient conditions on the residue field or value group
of a henselian valuation
to imply the existence of some nontrivial definable henselian valuation on that field.
The first one is about so-called {\em almost real closed} fields,
which is expressed as a condition on the residue field:

\begin{Theorem}[Delon--Farr\'e, {\cite[\S2]{DF}}]
Let $K$ be a field which is not real closed.
If $K$ admits a henselian valuation with real closed residue field,
then $K$ admits a nontrivial $\emptyset$-definable henselian valuation.
\footnote{Proof: There is a coarsest valuation $v_0$ on $K$ with real closed residue field \cite[Prop.~2.1]{DF},
whose valuation ring is the intersection of a family of valuation rings of valuations $v_{S}$, where $S$ runs over all 
finite sets of primes \cite[Proof of Lemma 2.7]{DF}.
By assumption, $v_0$ is nontrivial, hence also one of the $v_S$ must be nontrivial,
and it is $\emptyset$-definable by \cite[Prop.~2.6]{DF}.
}
\end{Theorem}
An analogous result holds when one replaces \emph{real closed} by \emph{separably closed}, see \cite{JK12}.
The second result in this direction uses a condition on the value group:

\begin{Theorem}[Jahnke--Koenigsmann, {\cite[Proposition 4.2]{JK16}}]
Let $K$ be a field.
If $K$ admits a henselian valuation with non-divisible value group,
then $K$ admits a nontrivial definable henselian valuation.
\end{Theorem}

Moreover, \cite[Theorems A and B]{JK16} give sufficient conditions on a henselian field $K$,
in terms of the value group and residue field of the canonical henselian valuation,
for $K$ to admit nontrivial definable resp.~$\emptyset$-definable henselian valuations.
In case the residue field of the canonical henselian valuation on $K$ has characteristic zero,
this is even a full characterization 
for the existence of nontrivial definable henselian valuations, cf.~\cite[Corollary 6.1]{JK16}.
As a concrete example they deduce that
every henselian field which is neither separably nor real closed and 
has finite transcendence degree over its prime field
admits a nontrivial $\emptyset$-definable henselian valuation.
In a different paper they 
give conditions on the absolute Galois group of a henselian field
that imply the existence of a nontrivial definable henselian valuation:

\begin{Theorem}[Jahnke--Koenigsmann, {\cite[Theorem 3.15]{JK12}}]
Let $K$ be a hensel\-ian field which is neither separably nor real closed.
If there exists a finite group $G$ such that no finite extension of $K$ has a Galois extension with Galois group $G$,
then $K$ admits a nontrivial $\emptyset$-definable henselian valuation.
\end{Theorem}

This result applies in particular to henselian valued {\em NIP} fields of positive characteristic that are not separably closed (see \cite[Corollary 3.18]{JK12}), which subsequently was generalized from NIP to $n$-NIP in \cite[Proposition 7.4]{Hempel}. For the definition of NIP and more details on NIP fields as well as on related model-theoretic concepts, see \cite{Sim15}.

Similarly, Johnson in \cite{Johnson} uses Theorem 
\ref{uniformp} to give a classification of {\em dp-minimal} fields which involves showing that 
every dp-minimal field which is neither algebraically closed or real closed admits a nontrivial definable henselian valuation.

Also \cite{Krupinski} gives some results towards the Shelah--Hasson conjecture: 
using \cite{Koe94}, he shows that every radically bounded field that admits a valuation with non-divisible value group
admits a nontrivial definable valuation (not necessarily henselian).
As an application, he deduces that every valuation on a 
field which is {\em superrosy} or {\em minimal} has divisible value group.

\section{Quantifier complexity and questions of uniformity}
\label{sec:3}

When defining sets in any structure, we are often interested in the quantifier complexity of the formulae involved and on determining
in which elementary classes a given definition uniformly yields the desired set. The classical results by J. Robinson (Theorem \ref{Robinson}) 
and Ax
(Theorem \ref{Ax}) discussed in
the first section give an $\exists$-$K$-formula and an $\exists\forall\exists\forall$-$\emptyset$-formula defining the valuation
ring $\mathcal{O}_v$ on any henselian valued field $(K,v)$ for which the value group is a $\mathbb{Z}$-group.
Moreover, one can check that the formulae worked out by Hong for $(K,v)$ henselian with a 
regular, non-discrete and non-divisible value groups (as
given in the proof of Theorem \ref{Hong} above) are in fact 
$\exists\forall$-$\emptyset$-formulae.

A class of henselian fields of particular interest is formed by the $p$-adics, their algebraic extensions and related fields
(e.g., $\mathbb{F}_p((t))$ and ultraproducts of the $p$-adics). Here,
we have the following results which were the starting point to many other recent works:
\begin{Theorem}[Cluckers--Derakhshan--Leenknegt--Macintyre, {\cite[Theorems 2, 5 and 6]{CDLM}}] \label{CDLM} \mbox{} 
\begin{enumerate} 
\item There is an $\exists\forall$-$\emptyset$-formula defining $\mathcal{O}_v$
in any henselian valued field $(K,v)$ with $Kv$ finite or pseudofinite. 
\item There is no $\forall$-$\emptyset$-formula nor an $\exists$-$\emptyset$-formula
which uniformly defines $\mathbb{Z}_p$ in $\mathbb{Q}_p$ for all primes $p$.
However, for any fixed finite extension $K$ of $\mathbb{Q}_p$, the
unique prolongation of $v_p$ to $K$ is both $\forall$-$\emptyset$-definable and 
$\exists$-$\emptyset$-definable.
\end{enumerate}
\end{Theorem}
Here, a field is called \emph{pseudofinite} if it is an infinite model of the common $\mathcal{L}_\mathrm{ring}$-theory
of all finite fields. 
In particular,
all pseudofinite fields are pseudo-algebraically closed but not separably closed. See \cite[\S20.10]{FJ} for more details
on pseudofinite fields.

In fact, the statement of Theorem \ref{CDLM}(1) given in \cite{CDLM} is stronger than the version given above, since they actually
prove $\exists$-$\emptyset$-definability in a modification of the Macintyre language $\mathcal{L}_\mathrm{Mac}$,
see the introduction of \cite{CDLM} for more details.
Note that if $(K,v)$ is henselian with finite (respectively pseudofinite) 
residue field and $L$ is a finite extension of $K$, 
then the residue field of the unique prolongation of $v$ to $L$ is again finite
(respectively pseudofinite). 
Thus, since for a fixed henselian valued field $(K,v)$ there is no 
$\forall$-$\emptyset$-formula nor an $\exists$-$\emptyset$-formula
which uniformly defines the unique prolongation of $v$ to $L$ for all finite extensions $L$ of $K$ (cf.~\cite[Theorem 4]{CDLM}),
the first statement of Theorem \ref{CDLM} cannot be improved to either 
a uniform $\forall$-$\emptyset$-definition or a uniform $\exists$-$\emptyset$-definition.

The positive characteristic analogue of Theorem \ref{CDLM}(2) was shown by Anscombe and Koenigsmann. Again, 
the definition cannot work uniformly for all power series fields over finite fields.
\begin{Theorem}[Anscombe--Koenigsmann, {\cite[Theorem 1.1]{AK}}] Let $q$ be a prime power. Then
$\mathbb{F}_q[[t]]$ is $\exists$-$\emptyset$-definable in $\mathbb{F}_q((t))$.
\end{Theorem}

Fehm generalizes the methods employed by Anscombe and Koenigsmann in \cite{F15} to work for all henselian valuations with
finite residue field. Moreover, he shows that although
uniformity for all primes is impossible to achieve, one can always find $\exists$-$\emptyset$-formulae which work
for large (infinite) families of finite residue fields:
\begin{Theorem}[Fehm, {\cite[Theorems 1.1 and 1.2]{F15}}] \label{Fehm} \mbox{} 
\begin{enumerate}
\item Let $(K,v)$ be henselian. If $Kv$ is finite or pseudo-algebraically closed and the algebraic part of $Kv$ is not algebraically
closed, then $\mathcal{O}_v$ is $\exists$-$\emptyset$-definable.
\item For every $\varepsilon >0$ there is an $\exists$-$\emptyset$-formula $\phi$ and a set $P$ of primes of
Dirichlet density at least $1- \varepsilon$ such that for any henselian $(K,v)$ with $|Kv| \in P$, the formula
$\phi$ defines $\mathcal{O}_v$ in $K$.
\end{enumerate}
\end{Theorem}
Note that the assumption on the algebraic part of $Kv$ being non-algebraically closed is necessary: The power series valuation
on $\mathbb{C}((t))$ is not $\exists$-$\emptyset$-definable (\cite[Appendix A]{AK}), 
and $\mathbb{C}$ is of course pseudo-algebraically closed.

All of the above results give explicit formulae. After these results had emerged, Prestel proved a Beth-like Characterization Theorem
in \cite{Prestel} which implies the existence of low-quantifier definitions (without a method to explicitly construct them).
Applying this, he shows the following (partial) improvement of Theorems \ref{JK12} and \ref{CDLM}:
\begin{Theorem}[Prestel, {\cite[Theorem 1]{Prestel}}]
There is an $\exists\forall$-$\emptyset$-formula defining uniformly the henselian valuations whose
residue field is either finite, pseudofinite or hilbertian.
\end{Theorem}

In \cite{FehmPrestel}, the authors apply Prestel's Characterization Theorem to obtain the existence of
uniform $\exists$-$\emptyset$- and
$\forall$-$\emptyset$-definitions for $\mathbb{Z}_p$ in $\mathbb{Q}_p$ and for $\mathbb{F}_p[[t]]$ in $\mathbb{F}_p((t))$
for odd primes $p$ in the Macintyre language $\mathcal{L}_\textrm{Mac}$. Moreover, they build on Hong's work (Theorem \ref{Hong}) to show
\begin{Theorem}[Fehm--Prestel, {\cite[Corollary 3.8]{FehmPrestel}}]
Let $(K,v)$ be a henselian valued field. 
If $vK$ is regular and non-divisible, then $v$ is $\exists\forall$-$\emptyset$-definable.
\end{Theorem}

Prestel's Characterization Theorem lead him to ask the question of 
whether whenever a henselian valuation is $\emptyset$-definable, 
it is
already definable with a formula of low quantifier complexity (that is with an $\exists\forall$-$\emptyset$-formula or an $\forall\exists$-$\emptyset$-formula). For canonical ($p$-)henselian valuations, this question is addressed
by Fehm and Jahnke in \cite{FJ15}. In the simpler case of the canonical $p$-henselian valuation and the setting of Theorem \ref{uniformp}, 
the canonical $p$-henselian
valuation is either $\forall\exists$-$\emptyset$-definable or $\exists\forall$-$\emptyset$-definable (depending on whether the residue
field is $p$-closed or not, see \cite[Propositions 3.6 and 3.7]{FJ15}).
For the canonical henselian valuation, the analogous result is only obtained if the absolute Galois group of the field is small 
(\cite[Theorem 1.1]{FJ15})

On the other hand, Halupczok and Jahnke construct an $\emptyset$-definable henselian valuation in \cite[Theorem 1.3]{HJ15} 
which is neither definable by an
$\exists\forall$-$\emptyset$-formula nor an $\forall\exists$-$\emptyset$-formula.

Very recently, Fehm and Anscombe gave a characterization of $\exists$-$\emptyset$-definable as well as $\forall$-$\emptyset$-definable henselian valuation rings
(\cite{AF2}). In particular, they show that the question of whether a given equicharacteristic henselian valuation 
on a field is $\exists$-$\emptyset$-definable 
(respectively $\forall$-$\emptyset$-definable)
depends only on the residue field of that valuation. The existential case of their theorem reads as follows.
\begin{Theorem}[Anscombe--Fehm, {\cite[Theorem 1.1]{AF2}}] Let $F$ be a field. Then the following are equivalent: \label{AF}
\begin{enumerate}
\item There is an $\exists$-$\emptyset$-formula that defines $\mathcal{O}_v$ in $K$ for \emph{some}
equicharacteristic henselian nontrivially valued field $(K,v)$ with residue field $F$.
\item There is an $\exists$-$\emptyset$-formula that defines $\mathcal{O}_v$ in $K$ for \emph{every}
henselian valued field $(K,v)$ with residue field elementarily equivalent to $F$.
\item There is no elementary extension $F^* \succ F$ with a nontrivial valuation $v$ on $F^*$ for which the
residue field $F^*v$ embeds into $F^*$.
\end{enumerate}
\end{Theorem} 
For the universal version of Theorem \ref{AF}, the quantifier $\exists$ is replaced by $\forall$ in conditions 
(1) and (2), and condition (3)
is replaced by
\begin{enumerate}
\item[(3)'] There is no elementary extension $F^* \succ F$ with a nontrivial henselian valuation $v$ on a subfield $E \subseteq F^*$
with $Ev \cong F^*$.
\end{enumerate} 
It is easy to see that both conditions (3) and (3)' hold for example in the setting of Theorem \ref{Fehm}(1).
Anscombe and Fehm also apply their results to show (\cite[Corollary 6.12]{AF2}) 
that for any field $F$, the valuation ring $F[[t]]$ is $\forall$-$F$-definable
on $F((t))$ if and only if $F$ is not a large field.
(\emph{Large} is a property
of fields which is a common generalization of henselian and pseudo-algebraically closed, 
see \cite[\S6.2]{AF2} for more details and references on large fields.)

\section*{Acknowledgements}
The authors would like to thank Martin Bays, Jochen Koenigsmann and Alexander Prestel for many helpful discussions and
their comments on this survey. 

\bibliographystyle{amsplain}

\end{document}